\documentstyle[11pt]{article}
\begin{document}
\newcommand{\p}{\parallel }
\makeatletter \makeatother
\newtheorem{th}{Theorem}[section]
\newtheorem{lem}{Lemma}[section]
\newtheorem{de}{Definition}[section]
\newtheorem{rem}{Remark}[section]
\newtheorem{cor}{Corollary}[section]
\renewcommand{\theequation}{\thesection.\arabic {equation}}

\title{{\bf Transgression and twisted anomaly cancellation formulas on odd dimensional manifolds}}

\author{ Yong Wang\\
 }

\date{}
\maketitle

\begin{abstract}~~ We compute the transgressed forms of some modularly invariant characteristic forms,
which are related to the twisted elliptic genera. We study the
modularity properties of these secondary characteristic forms and
relations among them.  We also get some twisted anomaly cancellation
formulas on some odd dimensional manifolds.\\

 \noindent{\bf Subj. Class.:}\quad Differential geometry; Algebraic topology\\
\noindent{\bf MSC:}\quad 58C20; 57R20; 53C80\\
 \noindent{\bf Keywords:}\quad
 Transgression; elliptic genera; cancellation formulas
\end{abstract}

\section{Introduction}
  \quad In 1983, the physicists Alvarez-Gaum\'{e} and Witten [AW]
  discovered the "miraculous cancellation" formula for gravitational
  anomaly which reveals a beautiful relation between the top
  components of the Hirzebruch $\widehat{L}$-form and
  $\widehat{A}$-form of a $12$-dimensional smooth Riemannian
  manifold. Kefeng Liu [Li] established higher dimensional "miraculous cancellation"
  formulas for $(8k+4)$-dimensional Riemannian manifolds by
  developing modular invariance properties of characteristic forms.
  These formulas could be used to deduce some divisibility results. In
  [HZ1], [HZ2], for each $(8k+4)$-dimensional smooth Riemannian
  manifold, a more general cancellation formula that involves a
  complex line bundle was established. This formula was applied to
  ${\rm spin}^c$ manifolds, then an analytic Ochanine congruence
  formula was derived. For $(8k+2)$ and $(8k+6)$-dimensional smooth Riemannian
  manifolds, F. Han and X. Huang [HH] obtained some cancellation formulas.\\
  \indent On the other hand, motivated by the Chern-Simons theory, in
  [CH], Qingtao Chen and Fei Han computed the transgressed forms of some modularly invariant characteristic forms,
which are related to the elliptic genera. They studied the
modularity properties of these secondary characteristic forms and
relations among them.  They also got an anomaly cancellation formula
for $11$-dimensional manifold. Thus a nature question is to get some
twisted modular forms by transgression and some twisted anomaly
cancellation formulas for odd dimensional manifolds. In this paper,
we compute the transgressed forms of some modularly invariant
characteristic forms, which are related to the "twisted" elliptic
genera. We study the modularity properties of these secondary
characteristic forms and relations among them.  We also get some
twisted anomaly cancellation formulas on some odd dimensional
manifolds. We hope that these new geometric invariants of
connections with modularity properties obtained here could be
applied somewhere.\\
\indent This paper is organized as follows: In Section 2, we review
some knowledge on characteristic forms and modular forms that we are
going to use. In Section 3, for $(4k-1)$ dimensional manifolds, we
apply the Chern-Simons transgression to characteristic forms with
modularity properties which are related to the "twisted" elliptic
genera and obtain some interesting secondary characteristic forms
with modularity properties. We also get two twisted cancellation
formulas for $11$-dimensional manifolds. In Section 4, for $(4k+1)$
dimensional manifolds, by transgression, we again obtain some
interesting secondary characteristic forms with modularity
properties. As a corollary, we get a twisted cancellation formula
for
$9$-dimensional manifolds.\\
\section{characteristic forms and modular forms}
 \quad The purpose of this section is to review the necessary knowledge on
characteristic forms and modular forms that we are going to use.\\

 \noindent {\bf  2.1 characteristic forms. }Let $M$ be a Riemannian manifold.
 Let $\nabla^{ TM}$ be the associated Levi-Civita connection on $TM$
 and $R^{TM}=(\nabla^{TM})^2$ be the curvature of $\nabla^{ TM}$.
 Let $\widehat{A}(TM,\nabla^{ TM})$ and $\widehat{L}(TM,\nabla^{ TM})$
 be the Hirzebruch characteristic forms defined respectively by (cf.
 [Z])
 $$\widehat{A}(TM,\nabla^{ TM})={\rm
 det}^{\frac{1}{2}}\left(\frac{\frac{\sqrt{-1}}{4\pi}R^{TM}}{{\rm
 sinh}(\frac{\sqrt{-1}}{4\pi}R^{TM})}\right),$$
 $$\widehat{L}(TM,\nabla^{ TM})={\rm
 det}^{\frac{1}{2}}\left(\frac{\frac{\sqrt{-1}}{2\pi}R^{TM}}{{\rm
 tanh}(\frac{\sqrt{-1}}{4\pi}R^{TM})}\right).\eqno(2.1)$$
   Let $E$, $F$ be two Hermitian vector bundles over $M$ carrying
   Hermitian connection $\nabla^E,\nabla^F$ respectively. Let
   $R^E=(\nabla^E)^2$ (resp. $R^F=(\nabla^F)^2$) be the curvature of
   $\nabla^E$ (resp. $\nabla^F$). If we set the formal difference
   $G=E-F$, then $G$ carries an induced Hermitian connection
   $\nabla^G$ in an obvious sense. We define the associated Chern
   character form as
   $${\rm ch}(G,\nabla^G)={\rm tr}\left[{\rm
   exp}(\frac{\sqrt{-1}}{2\pi}R^E)\right]-{\rm tr}\left[{\rm
   exp}(\frac{\sqrt{-1}}{2\pi}R^F)\right].\eqno(2.2)$$
   For any complex number $t$, let
   $$\wedge_t(E)={\bf C}|_M+tE+t^2\wedge^2(E)+\cdots,~S_t(E)={\bf
   C}|_M+tE+t^2S^2(E)+\cdots$$
   denote respectively the total exterior and symmetric powers of
   $E$, which live in $K(M)[[t]].$ The following relations between
   these operations hold,
   $$S_t(E)=\frac{1}{\wedge_{-t}(E)},~\wedge_t(E-F)=\frac{\wedge_t(E)}{\wedge_t(F)}.\eqno(2.3)$$
   Moreover, if $\{\omega_i\},\{\omega_j'\}$ are formal Chern roots
   for Hermitian vector bundles $E,F$ respectively, then
   $${\rm ch}(\wedge_t(E))=\prod_i(1+e^{\omega_i}t).\eqno(2.4)$$
   Then we have the following formulas for Chern character forms,
   $${\rm ch}(S_t(E))=\frac{1}{\prod_i(1-e^{\omega_i}t)},~
{\rm
ch}(\wedge_t(E-F))=\frac{\prod_i(1+e^{\omega_i}t)}{\prod_j(1+e^{\omega'_j}t)}.\eqno(2.5)$$
\indent If $W$ is a real Euclidean vector bundle over $M$ carrying a
Euclidean connection $\nabla^W$, then its complexification $W_{\bf
C}=W\otimes {\bf C}$ is a complex vector bundle over $M$ carrying a
canonical induced Hermitian metric from that of $W$, as well as a
Hermitian connection $\nabla^{W_{\bf C}}$ induced from $\nabla^W$.
If $E$ is a vector bundle (complex or real) over $M$, set
$\widetilde{E}=E-{\rm dim}E$ in $K(M)$ or $KO(M)$.\\

\noindent{\bf 2.2 Some properties about the Jacobi theta functions
and modular forms}\\
   \indent We first recall the four Jacobi theta functions are
   defined as follows( cf. [Ch]):
   $$\theta(v,\tau)=2q^{\frac{1}{8}}{\rm sin}(\pi
   v)\prod_{j=1}^{\infty}[(1-q^j)(1-e^{2\pi\sqrt{-1}v}q^j)(1-e^{-2\pi\sqrt{-1}v}q^j)],\eqno(2.6)$$
$$\theta_1(v,\tau)=2q^{\frac{1}{8}}{\rm cos}(\pi
   v)\prod_{j=1}^{\infty}[(1-q^j)(1+e^{2\pi\sqrt{-1}v}q^j)(1+e^{-2\pi\sqrt{-1}v}q^j)],\eqno(2.7)$$
$$\theta_2(v,\tau)=\prod_{j=1}^{\infty}[(1-q^j)(1-e^{2\pi\sqrt{-1}v}q^{j-\frac{1}{2}})
(1-e^{-2\pi\sqrt{-1}v}q^{j-\frac{1}{2}})],\eqno(2.8)$$
$$\theta_3(v,\tau)=\prod_{j=1}^{\infty}[(1-q^j)(1+e^{2\pi\sqrt{-1}v}q^{j-\frac{1}{2}})
(1+e^{-2\pi\sqrt{-1}v}q^{j-\frac{1}{2}})],\eqno(2.9)$$ \noindent
where $q=e^{2\pi\sqrt{-1}\tau}$ with $\tau\in\textbf{H}$, the upper
half complex plane. Let
$$\theta'(0,\tau)=\frac{\partial\theta(v,\tau)}{\partial
v}|_{v=0}.\eqno(2.10)$$ \noindent Then the following Jacobi identity
(cf. [Ch]) holds,
$$\theta'(0,\tau)=\pi\theta_1(0,\tau)\theta_2(0,\tau)\theta_3(0,\tau).\eqno(2.11)$$
\noindent Denote $SL_2({\bf Z})=\left\{\left(\begin{array}{cc}
\ a & b  \\
 c  & d
\end{array}\right)\mid a,b,c,d \in {\bf Z},~ad-bc=1\right\}$ the
modular group. Let $S=\left(\begin{array}{cc}
\ 0 & -1  \\
 1  & 0
\end{array}\right),~T=\left(\begin{array}{cc}
\ 1 &  1 \\
 0  & 1
\end{array}\right)$ be the two generators of $SL_2(\bf{Z})$. They
act on $\textbf{H}$ by $S\tau=-\frac{1}{\tau},~T\tau=\tau+1$. One
has the following transformation laws of theta functions under the
actions of $S$ and $T$ (cf. [Ch]):
$$\theta(v,\tau+1)=e^{\frac{\pi\sqrt{-1}}{4}}\theta(v,\tau),~~\theta(v,-\frac{1}{\tau})
=\frac{1}{\sqrt{-1}}\left(\frac{\tau}{\sqrt{-1}}\right)^{\frac{1}{2}}e^{\pi\sqrt{-1}\tau
v^2}\theta(\tau v,\tau);\eqno(2.12)$$
$$\theta_1(v,\tau+1)=e^{\frac{\pi\sqrt{-1}}{4}}\theta_1(v,\tau),~~\theta_1(v,-\frac{1}{\tau})
=\left(\frac{\tau}{\sqrt{-1}}\right)^{\frac{1}{2}}e^{\pi\sqrt{-1}\tau
v^2}\theta_2(\tau v,\tau);\eqno(2.13)$$
$$\theta_2(v,\tau+1)=\theta_3(v,\tau),~~\theta_2(v,-\frac{1}{\tau})
=\left(\frac{\tau}{\sqrt{-1}}\right)^{\frac{1}{2}}e^{\pi\sqrt{-1}\tau
v^2}\theta_1(\tau v,\tau);\eqno(2.14)$$
$$\theta_3(v,\tau+1)=\theta_2(v,\tau),~~\theta_3(v,-\frac{1}{\tau})
=\left(\frac{\tau}{\sqrt{-1}}\right)^{\frac{1}{2}}e^{\pi\sqrt{-1}\tau
v^2}\theta_3(\tau v,\tau).\eqno(2.15)$$ \noindent Differentiating
the above transformation formulas, we get that
$$\theta'(v,\tau+1)=e^{\frac{\pi\sqrt{-1}}{4}}\theta'(v,\tau),$$
$$\theta'(v,-\frac{1}{\tau})=\frac{1}{\sqrt{-1}}\left(\frac{\tau}{\sqrt{-1}}\right)^{\frac{1}{2}}e^{\pi\sqrt{-1}\tau
v^2}(2\pi\sqrt{-1}\tau v\theta(\tau v,\tau)+\tau\theta'(\tau
v,\tau));$$
$$\theta'_1(v,\tau+1)=e^{\frac{\pi\sqrt{-1}}{4}}\theta_1'(v,\tau),$$
$$\theta'_1(v,-\frac{1}{\tau})=\left(\frac{\tau}{\sqrt{-1}}\right)^{\frac{1}{2}}e^{\pi\sqrt{-1}\tau
v^2}(2\pi\sqrt{-1}\tau v\theta_2(\tau v,\tau)+\tau\theta'_2(\tau
v,\tau));$$
$$\theta'_2(v,\tau+1)=\theta_3'(v,\tau),$$
$$\theta'_2(v,-\frac{1}{\tau})=\left(\frac{\tau}{\sqrt{-1}}\right)^{\frac{1}{2}}e^{\pi\sqrt{-1}\tau
v^2}(2\pi\sqrt{-1}\tau v\theta_1(\tau v,\tau)+\tau\theta'_1(\tau
v,\tau));$$
$$\theta'_3(v,\tau+1)=\theta_2'(v,\tau),$$
$$\theta'_3(v,-\frac{1}{\tau})=\left(\frac{\tau}{\sqrt{-1}}\right)^{\frac{1}{2}}e^{\pi\sqrt{-1}\tau
v^2}(2\pi\sqrt{-1}\tau v\theta_3(\tau v,\tau)+\tau\theta'_3(\tau
v,\tau))\eqno(2.16)$$
 \noindent Therefore
 $$\theta'(0,-\frac{1}{\tau})=\frac{1}{\sqrt{-1}}\left(\frac{\tau}{\sqrt{-1}}\right)^{\frac{1}{2}}
\tau\theta'(0,\tau).\eqno(2.17)$$
 \noindent {\bf Definition 2.1} A modular form over $\Gamma$, a
 subgroup of $SL_2({\bf Z})$, is a holomorphic function $f(\tau)$ on
 $\textbf{H}$ such that
 $$f(g\tau):=f\left(\frac{a\tau+b}{c\tau+d}\right)=\chi(g)(c\tau+d)^kf(\tau),
 ~~\forall g=\left(\begin{array}{cc}
\ a & b  \\
 c & d
\end{array}\right)\in\Gamma,\eqno(2.18)$$
\noindent where $\chi:\Gamma\rightarrow {\bf C}^{\star}$ is a
character of $\Gamma$. $k$ is called the weight of $f$.\\
Let $$\Gamma_0(2)=\left\{\left(\begin{array}{cc}
\ a & b  \\
 c  & d
\end{array}\right)\in SL_2({\bf Z})\mid c\equiv 0~({\rm
mod}~2)\right\},$$
$$\Gamma^0(2)=\left\{\left(\begin{array}{cc}
\ a & b  \\
 c  & d
\end{array}\right)\in SL_2({\bf Z})\mid b\equiv 0~({\rm
mod}~2)\right\},$$
$$\Gamma_\theta=\left\{\left(\begin{array}{cc}
\ a & b  \\
 c  & d
\end{array}\right)\in SL_2({\bf Z})\mid
\left(\begin{array}{cc}
\ a & b  \\
 c  & d
\end{array}\right)
\equiv \left(\begin{array}{cc}
\ 1 & 0  \\
 0  & 1
\end{array}\right)
{\rm or} \left(\begin{array}{cc}
\ 0 & 1  \\
 1  & 0
\end{array}\right)
 ~({\rm
mod}~2)\right\}$$ be the three modular subgroups of $SL_2({\bf Z})$.
It is known that the generators of $\Gamma_0(2)$ are $T,~ST^2ST$,
the generators of $\Gamma^0(2)$ are $STS,~T^2STS$ and the generators
of $\Gamma_\theta$ are $S,~T^2$ (cf.[Ch]).\\
\indent If $\Gamma$ is a modular subgroup, let ${\mathcal{M}}_{{\bf
R}}(\Gamma)$ denote the ring of modular forms over $\Gamma$ with
real Fourier coefficients. Writing $\theta_j=\theta_j(0,\tau),~1\leq
j\leq 3,$ we introduce six explicit modular forms (cf. [Li]),
$$\delta_1(\tau)=\frac{1}{8}(\theta_2^4+\theta_3^4),~~\varepsilon_1(\tau)=\frac{1}{16}\theta_2^4\theta_3^4,$$
$$\delta_2(\tau)=-\frac{1}{8}(\theta_1^4+\theta_3^4),~~\varepsilon_2(\tau)=\frac{1}{16}\theta_1^4\theta_3^4,$$
$$\delta_3(\tau)=\frac{1}{8}(\theta_1^4-\theta_2^4),~~\varepsilon_3(\tau)=-\frac{1}{16}\theta_1^4\theta_2^4.$$
\noindent They have the following Fourier expansions in
$q^{\frac{1}{2}}$:
$$\delta_1(\tau)=\frac{1}{4}+6q+\cdots,~~\varepsilon_1(\tau)=\frac{1}{16}-q+\cdots,$$
$$\delta_2(\tau)=-\frac{1}{8}-3q^{\frac{1}{2}}+\cdots,~~\varepsilon_2(\tau)=q^{\frac{1}{2}}+\cdots,$$
$$\delta_3(\tau)=-\frac{1}{8}+3q^{\frac{1}{2}}+\cdots,~~\varepsilon_3(\tau)=-q^{\frac{1}{2}}+\cdots,$$
\noindent where the $"\cdots"$ terms are the higher degree terms,
all of which have integral coefficients. They also satisfy the
transformation laws,
$$\delta_2(-\frac{1}{\tau})=\tau^2\delta_1(\tau),~~~~~~\varepsilon_2(-\frac{1}{\tau})
=\tau^4\varepsilon_1(\tau),\eqno(2.19)$$
$$\delta_2(\tau+1)=\delta_3(\tau),~~~~~~\varepsilon_2(\tau+1)=\varepsilon_3(\tau).\eqno(2.20)$$
\noindent {\bf Lemma 2.2} ([Li]) {\it $\delta_1(\tau)$ (resp.
$\varepsilon_1(\tau)$) is a modular form of weight $2$ (resp. $4$)
over $\Gamma_0(2)$, $\delta_2(\tau)$ (resp. $\varepsilon_2(\tau)$)
is a modular form of weight $2$ (resp. $4$) over $\Gamma^0(2)$,
while  $\delta_3(\tau)$ (resp. $\varepsilon_3(\tau)$) is a modular
form of weight $2$ (resp. $4$) over $\Gamma_\theta(2)$ and moreover
${\mathcal{M}}_{{\bf R}}(\Gamma^0(2))={\bf
R}[\delta_2(\tau),\varepsilon_2(\tau)]$.}

\section {Transgressed forms and modularities on $4k-1$ dimensional
manifolds}

   \quad Let $M$ be a $4k-1$ dimensional Riemannian
manifold and
 $\xi$ be a rank two real oriented Euclidean vector
   bundle over $M$ carrying with a Euclidean connection
   $\nabla^\xi$. Set
   $$\Theta_1(T_{C}M,\xi_C)=\bigotimes _{n=1}^{\infty}S_{q^n}(\widetilde{T_CM})\otimes
\bigotimes
_{m=1}^{\infty}\wedge_{q^m}(\widetilde{T_CM}-2\widetilde{\xi_C})\otimes
\bigotimes _{r=1}^{\infty}\wedge
_{q^{r-\frac{1}{2}}}(\widetilde{\xi_C})\otimes\bigotimes
_{s=1}^{\infty}\wedge _{-q^{s-\frac{1}{2}}}(\widetilde{\xi_C}),$$
$$\Theta_2(T_{C}M,\xi_C)=\bigotimes _{n=1}^{\infty}S_{q^n}(\widetilde{T_CM})\otimes
\bigotimes
_{m=1}^{\infty}\wedge_{-q^{m-\frac{1}{2}}}(\widetilde{T_CM}-2\widetilde{\xi_C})\otimes
\bigotimes _{r=1}^{\infty}\wedge
_{q^{r-\frac{1}{2}}}(\widetilde{\xi_C})\otimes\bigotimes
_{s=1}^{\infty}\wedge _{q^{s}}(\widetilde{\xi_C}),$$
$$\Theta_3(T_{C}M,\xi_C)=\bigotimes _{n=1}^{\infty}S_{q^n}(\widetilde{T_CM})\otimes
\bigotimes
_{m=1}^{\infty}\wedge_{q^{m-\frac{1}{2}}}(\widetilde{T_CM}-2\widetilde{\xi_C})\otimes
\bigotimes _{r=1}^{\infty}\wedge
_{q^{r}}(\widetilde{\xi_C})\otimes\bigotimes _{s=1}^{\infty}\wedge
_{-q^{s-\frac{1}{2}}}(\widetilde{\xi_C}).\eqno(3.1)$$ Let
$c=e(\xi,\nabla^{\xi})$ be the Euler form of $\xi$ canonically
associated to $\nabla^\xi$. Set
$$\Phi_L(\nabla^{TM},\nabla^{\xi},\tau)=\frac{\widehat{L}(TM,\nabla^{TM})}{{\rm
cosh}^2(\frac{c}{2})}{\rm
ch}(\Theta_1(T_{C}M,\xi_C),\nabla^{\Theta_1(T_{C}M,\xi_C)}),$$
$$\Phi_W(\nabla^{TM},\nabla^{\xi},\tau)={\widehat{A}(TM,\nabla^{TM})}{\rm
cosh}(\frac{c}{2}){\rm
ch}(\Theta_2(T_{C}M,\xi_C),\nabla^{\Theta_2(T_{C}M,\xi_C)}),$$
$$\Phi_W'(\nabla^{TM},\nabla^{\xi},\tau)={\widehat{A}(TM,\nabla^{TM})}{\rm
cosh}(\frac{c}{2}){\rm
ch}(\Theta_3(T_{C}M,\xi_C),\nabla^{\Theta_3(T_{C}M,\xi_C)}).\eqno(3.2)$$
Let $\{\pm2\pi\sqrt{-1}x_j| ~1\leq j\leq 2k-1\}$ and
$\{\pm2\pi\sqrt{-1}u\}$ be the Chern roots of $T_CM$ and $\xi_C$
respectively and $c=2\pi\sqrt{-1}u.$ Through direct computations, we
get (cf. [HZ2])
$$\Phi_L(\nabla^{TM},\nabla^{\xi},\tau)=\sqrt{2}^{4k-1}\left\{
\left(\prod_{j=1}^{2k-1}x_j\frac{\theta'(0,\tau)}{\theta(x_j,\tau)}\frac
{\theta_1(x_j,\tau)}{\theta_1(0,\tau)}\right)\frac{\theta_1^2(0,\tau)}
{\theta_1^2(u,\tau)}\frac{\theta_3(u,\tau)}{\theta_3(0,\tau)}
\frac{\theta_2(u,\tau)}{\theta_2(0,\tau)}\right\};\eqno(3.3)$$
$$\Phi_W(\nabla^{TM},\nabla^{\xi},\tau)=\left(
\prod_{j=1}^{2k-1}x_j\frac{\theta'(0,\tau)}{\theta(x_j,\tau)}\frac
{\theta_2(x_j,\tau)}{\theta_2(0,\tau)}\right)\frac{\theta_2^2(0,\tau)}
{\theta_2^2(u,\tau)}\frac{\theta_3(u,\tau)}{\theta_3(0,\tau)}
\frac{\theta_1(u,\tau)}{\theta_1(0,\tau)};\eqno(3.4)$$
$$\Phi_W'(\nabla^{TM},\nabla^{\xi},\tau)=\left(
\prod_{j=1}^{2k-1}x_j\frac{\theta'(0,\tau)}{\theta(x_j,\tau)}\frac
{\theta_3(x_j,\tau)}{\theta_3(0,\tau)}\right)\frac{\theta_3^2(0,\tau)}
{\theta_3^2(u,\tau)}\frac{\theta_1(u,\tau)}{\theta_1(0,\tau)}
\frac{\theta_2(u,\tau)}{\theta_2(0,\tau)}.\eqno(3.5)$$ Consider the
following function defined on ${\bf C}\times {\bf H}$,
$$f_{\Phi_L}(z,\tau)=z\frac{\theta'(0,\tau)}{\theta(z,\tau)}\frac
{\theta_1(z,\tau)}{\theta_1(0,\tau)},$$
$$f_{\Phi_W}(z,\tau)=z\frac{\theta'(0,\tau)}{\theta(z,\tau)}\frac
{\theta_2(z,\tau)}{\theta_2(0,\tau)},$$
$$f_{\Phi_W'}(z,\tau)=z\frac{\theta'(0,\tau)}{\theta(z,\tau)}\frac
{\theta_3(z,\tau)}{\theta_3(0,\tau)}.$$ Applying the Chern-Weil
theory, we can express $\Phi_L,~\Phi_W, \Phi_W'$ as follows:
$$\Phi_L(\nabla^{TM},\nabla^{\xi},\tau)=\sqrt{2}^{4k-1}
{\rm
det}^{\frac{1}{2}}\left(f_{\Phi_L}(\frac{R^{TM}}{4\pi^2},\tau)\right)
{\rm det}^{\frac{1}{2}}\left(\frac{\theta_1^2(0,\tau)}
{\theta_1^2(\frac{R^{\xi}}{4\pi^2},\tau)}\frac{\theta_3(\frac{R^{\xi}}{4\pi^2},\tau)}{\theta_3(0,\tau)}
\frac{\theta_2(\frac{R^{\xi}}{4\pi^2},\tau)}{\theta_2(0,\tau)}\right);\eqno(3.6)$$
$$\Phi_W(\nabla^{TM},\nabla^{\xi},\tau)=
{\rm
det}^{\frac{1}{2}}\left(f_{\Phi_W}(\frac{R^{TM}}{4\pi^2},\tau)\right)
{\rm det}^{\frac{1}{2}}\left(\frac{\theta_2^2(0,\tau)}
{\theta_2^2(\frac{R^{\xi}}{4\pi^2},\tau)}\frac{\theta_3(\frac{R^{\xi}}{4\pi^2},\tau)}{\theta_3(0,\tau)}
\frac{\theta_1(\frac{R^{\xi}}{4\pi^2},\tau)}{\theta_1(0,\tau)}\right);\eqno(3.7)$$
$$\Phi_W'(\nabla^{TM},\nabla^{\xi},\tau)=
{\rm
det}^{\frac{1}{2}}\left(f_{\Phi_W'}(\frac{R^{TM}}{4\pi^2},\tau)\right)
{\rm det}^{\frac{1}{2}}\left(\frac{\theta_3^2(0,\tau)}
{\theta_3^2(\frac{R^{\xi}}{4\pi^2},\tau)}\frac{\theta_1(\frac{R^{\xi}}{4\pi^2},\tau)}{\theta_1(0,\tau)}
\frac{\theta_2(\frac{R^{\xi}}{4\pi^2},\tau)}{\theta_2(0,\tau)}\right).\eqno(3.8)$$
Let $E$ be a vector bundle and $f$ be a power series with constant
term $1$. Let $\nabla_t^E$ be deformed connection given by
$\nabla_t^E=(1-t)\nabla_0^E+t\nabla_1^E$ and $R^E_t, ~t\in [0,1],$
denote the curvature of $\nabla_t^E$. $f'(t)$ is the power series
obtained from the derivative of $f(x)$ with respect to $x$. $\omega$
is a closed form. Recall the trivial modification of Theorem 2.2 in
[CH],\\

 \noindent {\bf Lemma 3.1}( [CH])
$${\rm det}^{\frac{1}{2}}(f(R^E_1))\omega-{\rm det}^{\frac{1}{2}}(f(R^E_0))\omega
=d\int_0^1\frac{1}{2}{\rm det}^{\frac{1}{2}}(f(R^E_t))\omega{\rm
tr}\left[\frac{d\nabla_t^E}{dt}\frac{f'(R^E_t)}{f(R^E_t)}\right]dt.\eqno(3.9)$$
Now we let $E=TM$ and $A=\nabla_1^{TM}-\nabla_0^{TM}$, then by Lemma
3.1, we have \\
$\Phi_L(\nabla_1^{TM},\nabla^{\xi},\tau)
-\Phi_L(\nabla_0^{TM},\nabla^{\xi},\tau)$ $$=
\frac{1}{8\pi^2}d\int_0^1
\Phi_L(\nabla_t^{TM},\nabla^{\xi},\tau){\rm tr}\left[A \left(
\frac{1}{\frac{R_t^{TM}}{4\pi^2}}-\frac{\theta'(\frac{R_t^{TM}}{4\pi^2},\tau)}
{\theta(\frac{R_t^{TM}}{4\pi^2},\tau)}+\frac{\theta_1'(\frac{R_t^{TM}}{4\pi^2},\tau)}
{\theta_1(\frac{R_t^{TM}}{4\pi^2},\tau)}\right)\right]dt.\eqno(3.10)$$
We define\\
$CS\Phi_L(\nabla_0^{TM},\nabla_1^{TM},\nabla^{\xi},\tau)$
$$:= \frac{\sqrt{2}}{8\pi^2}\int_0^1
\Phi_L(\nabla_t^{TM},\nabla^{\xi},\tau){\rm tr}\left[A \left(
\frac{1}{\frac{R_t^{TM}}{4\pi^2}}-\frac{\theta'(\frac{R_t^{TM}}{4\pi^2},\tau)}
{\theta(\frac{R_t^{TM}}{4\pi^2},\tau)}+\frac{\theta_1'(\frac{R_t^{TM}}{4\pi^2},\tau)}
{\theta_1(\frac{R_t^{TM}}{4\pi^2},\tau)}\right)\right]dt.\eqno(3.11)$$
which is in $\Omega^{\rm odd}(M,{\bf C})[[q^{\frac{1}{2}}]].$ Since
$M$ is $4k-1$ dimensional,
$\{CS\Phi_L(\nabla_0^{TM},\nabla_1^{TM},\nabla^{\xi},\tau)\}^{(4k-1)}$
represents an element in $H^{4k-1}(M,{\bf C})[[q^{\frac{1}{2}}]]$.
Similarly, we can compute the transgressed forms for
$\Phi_W,~\Phi_W'$ respectively and define\\
$CS\Phi_W(\nabla_0^{TM},\nabla_1^{TM},\nabla^{\xi},\tau)$
$$:= \frac{1}{8\pi^2}\int_0^1
\Phi_W(\nabla_t^{TM},\nabla^{\xi},\tau){\rm tr}\left[A \left(
\frac{1}{\frac{R_t^{TM}}{4\pi^2}}-\frac{\theta'(\frac{R_t^{TM}}{4\pi^2},\tau)}
{\theta(\frac{R_t^{TM}}{4\pi^2},\tau)}+\frac{\theta_2'(\frac{R_t^{TM}}{4\pi^2},\tau)}
{\theta_2(\frac{R_t^{TM}}{4\pi^2},\tau)}\right)\right]dt;\eqno(3.12)$$
$CS\Phi_W'(\nabla_0^{TM},\nabla_1^{TM},\nabla^{\xi},\tau)$
$$:= \frac{1}{8\pi^2}\int_0^1
\Phi_W'(\nabla_t^{TM},\nabla^{\xi},\tau){\rm tr}\left[A \left(
\frac{1}{\frac{R_t^{TM}}{4\pi^2}}-\frac{\theta'(\frac{R_t^{TM}}{4\pi^2},\tau)}
{\theta(\frac{R_t^{TM}}{4\pi^2},\tau)}+\frac{\theta_3'(\frac{R_t^{TM}}{4\pi^2},\tau)}
{\theta_3(\frac{R_t^{TM}}{4\pi^2},\tau)}\right)\right]dt,\eqno(3.13)$$
which also lie in $\Omega^{\rm odd}(M,{\bf C})[[q^{\frac{1}{2}}]]$
and their top components represent elements in $H^{4k-1}(M,{\bf
C})[[q^{\frac{1}{2}}]]$. As pointed in [CH], the equality (3.10) and
the modular invariance properties of
$\Phi_L(\nabla_0^{TM},\nabla^{\xi},\tau)$ and
$\Phi_L(\nabla_1^{TM},\nabla^{\xi},\tau)$ are not enough to
guarantee that
$CS\Phi_L(\nabla_0^{TM},\nabla_1^{TM},\nabla^{\xi},\tau)$ is a
modular form. However we have the following results.\\

\noindent {\bf Theorem 3.2} {\it Let $M$ be a $4k-1$ dimensional
manifold and  $\nabla_0^{TM},~\nabla_1^{TM}$ be two connections on
$TM$ and $\xi$ be a two dimensional oriented Euclidean real vector
bundle with a Euclidean connection $\nabla^{\xi}$, then we have\\
\noindent 1)
$\{CS\Phi_L(\nabla_0^{TM},\nabla_1^{TM},\nabla^{\xi},\tau)\}^{(4k-1)}$
is a modular form of weight $2k$ over $\Gamma_0(2)$;\\
$\{CS\Phi_W(\nabla_0^{TM},\nabla_1^{TM},\nabla^{\xi},\tau)\}^{(4k-1)}$
is a modular form of weight $2k$ over $\Gamma^0(2);$\\
$\{CS\Phi_W'(\nabla_0^{TM},\nabla_1^{TM},\nabla^{\xi},\tau)\}^{(4k-1)}$
is a modular form of weight $2k$ over $\Gamma_\theta(2).$\\
2) The following equalities hold,}
$$\{CS\Phi_L(\nabla_0^{TM},\nabla_1^{TM},\nabla^{\xi},-\frac{1}{\tau})\}^{(4k-1)}
=(2\tau)^{2k}\{CS\Phi_W(\nabla_0^{TM},\nabla_1^{TM},\nabla^{\xi},\tau)\}^{(4k-1)},$$
$$CS\Phi_W(\nabla_0^{TM},\nabla_1^{TM},\nabla^{\xi},\tau+1)
=CS\Phi_W'(\nabla_0^{TM},\nabla_1^{TM},\nabla^{\xi},\tau).$$
\noindent{\bf Proof.} By (2.12)-(2.17), we have
$$z\frac{\theta'(0,-\frac{1}{\tau})}{\theta(z,-\frac{1}{\tau})}\frac
{\theta_1(z,-\frac{1}{\tau})}{\theta_1(0,-\frac{1}{\tau})} =(\tau
z)\frac{\theta'(0,{\tau})}{\theta(\tau z,{\tau})}\frac
{\theta_2(\tau z,{\tau})}{\theta_2(0,{\tau})};\eqno(3.14)$$
$$\frac{1}{z}-\frac{\theta'(z,-\frac{1}{\tau})}
{\theta(z,-\frac{1}{\tau})}+\frac{\theta_1'(z,-\frac{1}{\tau})}{\theta_1(z,-\frac{1}{\tau})}
=\tau\left(\frac{1}{\tau z}-\frac{\theta'(\tau z,{\tau})}
{\theta(\tau z,{\tau})}+\frac{\theta_2'(\tau
z,{\tau})}{\theta_2(\tau z,{\tau})} \right);\eqno(3.15)$$
$$\frac{\theta_1^2(0,-\frac{1}{\tau})}
{\theta_1^2(u,-\frac{1}{\tau})}\frac{\theta_3(u,-\frac{1}{\tau})}{\theta_3(0,-\frac{1}{\tau})}
\frac{\theta_2(u,-\frac{1}{\tau})}{\theta_2(0,-\frac{1}{\tau})}
=\frac{\theta_2^2(0,{\tau})} {\theta_2^2(\tau
u,{\tau})}\frac{\theta_3(\tau u,{\tau})}{\theta_3(0,{\tau})}
\frac{\theta_1(u\tau,{\tau})}{\theta_1(0,{\tau})}.\eqno(3.16)$$ Note
that we only take $(4k-1)$-component, so by (3.6)-(3.8),(3.11),
(3.12), (3.14)-(3.16), we can get
$$\{CS\Phi_L(\nabla_0^{TM},\nabla_1^{TM},\nabla^{\xi},-\frac{1}{\tau})\}^{(4k-1)}
=(2\tau)^{2k}\{CS\Phi_W(\nabla_0^{TM},\nabla_1^{TM},\nabla^{\xi},\tau)\}^{(4k-1)},\eqno(3.17)$$
Similarly we can show that
$$CS\Phi_L(\nabla_0^{TM},\nabla_1^{TM},\nabla^{\xi},{\tau}+1)
=CS\Phi_L(\nabla_0^{TM},\nabla_1^{TM},\nabla^{\xi},\tau),$$
$$\{CS\Phi_W(\nabla_0^{TM},\nabla_1^{TM},\nabla^{\xi},-\frac{1}{\tau})\}^{(4k-1)}
=(\frac{\tau}{2})^{2k}\{CS\Phi_L(\nabla_0^{TM},\nabla_1^{TM},\nabla^{\xi},\tau)\}^{(4k-1)},$$
$$CS\Phi_W(\nabla_0^{TM},\nabla_1^{TM},\nabla^{\xi},{\tau}+1)
=CS\Phi_W'(\nabla_0^{TM},\nabla_1^{TM},\nabla^{\xi},\tau),$$
$$\{CS\Phi_W'(\nabla_0^{TM},\nabla_1^{TM},\nabla^{\xi},-\frac{1}{\tau})\}^{(4k-1)}
=(\tau)^{2k}\{CS\Phi_W'(\nabla_0^{TM},\nabla_1^{TM},\nabla^{\xi},\tau)\}^{(4k-1)},$$
$$CS\Phi_W'(\nabla_0^{TM},\nabla_1^{TM},\nabla^{\xi},{\tau}+1)
=CS\Phi_W(\nabla_0^{TM},\nabla_1^{TM},\nabla^{\xi},\tau).\eqno(3.18)$$
From (3.17) and (3.18), we can get
$\{CS\Phi_L(\nabla_0^{TM},\nabla_1^{TM},\nabla^{\xi},{\tau})\}^{(4k-1)}$
is a modular form of weight $2k$ over $\Gamma_0(2).$ Similarly we
can prove that
$\{CS\Phi_W(\nabla_0^{TM},\nabla_1^{TM},\nabla^{\xi},\tau)\}^{(4k-1)}$
is a modular form of weight $2k$ over $\Gamma^0(2)$ and
$\{CS\Phi_W'(\nabla_0^{TM},\nabla_1^{TM},\nabla^{\xi},\tau)\}^{(4k-1)}$
is a modular form of weight $2k$ over $\Gamma_\theta(2).$ $\Box$\\
 \indent Let $M$ be a compact oriented smooth $3$-dimensional manifold, then
 our transgressed forms are same as transgressed forms in the untwisted
case which have been computed in [CH]. From Theorem 3.2, we can
imply some twisted cancellation formulas for odd dimensional
manifolds. For example, let $M$ be $11$ dimensional and $k=3$. We
have that
$\{CS\Phi_L(\nabla_0^{TM},\nabla_1^{TM},\nabla^{\xi},\tau)\}^{(11)}$
is a modular form of weight $6$ over $\Gamma_0(2),$
$\{CS\Phi_W(\nabla_0^{TM},\nabla_1^{TM},\nabla^{\xi},\tau)\}^{(11)}$
is a modular form of weight $6$ over $\Gamma^0(2)$ and
$$\{CS\Phi_L(\nabla_0^{TM},\nabla_1^{TM},\nabla^{\xi},-\frac{1}{\tau})\}^{(11)}
=(2\tau)^{6}\{CS\Phi_W(\nabla_0^{TM},\nabla_1^{TM},\nabla^{\xi},\tau)\}^{(11)}.$$
By Lemma 2.2, we have
$$\{CS\Phi_W(\nabla_0^{TM},\nabla_1^{TM},\nabla^{\xi},\tau)\}^{(11)}
=z_0(8\delta_2)^3+z_1(8\delta_2)\varepsilon_2,\eqno(3.19)$$ and by
(2.19) and Theorem 3.2,
$$\{CS\Phi_L(\nabla_0^{TM},\nabla_1^{TM},\nabla^{\xi},\tau)\}^{(11)}
=2^6[z_0(8\delta_1)^3+z_1(8\delta_1)\varepsilon_1].\eqno(3.20)$$ By
comparing the $q^{\frac{1}{2}}$-expansion coefficients in (3.19), we
get
$$z_0=-\left\{\int_0^1\widehat{A}(TM,\nabla^{TM}_t){\rm
cosh}(\frac{c}{2}){\rm
tr}\left[A\left(\frac{1}{2R_t^{TM}}-\frac{1}{8\pi{\rm
tan}{\frac{R^{TM}_t}{4\pi}}}\right)\right]dt\right\}^{(11)},\eqno(3.21)$$
$$z_1=\left\{\int_0^1\widehat{A}(TM,\nabla^{TM}_t){\rm
cosh}(\frac{c}{2})\left({\rm
ch}(T_CM,\nabla^{T_CM}_t)-3(e^c+e^{-c}-2)\right)\right.$$
$$\left.\times{\rm tr}\left[A\left(\frac{1}{2R_t^{TM}}-\frac{1}{8\pi{\rm
tan}{\frac{R^{TM}_t}{4\pi}}}\right)\right]dt\right.+$$
$$\left.\int_0^1\widehat{A}(TM,\nabla^{TM}_t){\rm
cosh}(\frac{c}{2}){\rm tr}\left[A\left(-\frac{1}{2\pi}{\rm
sin}\frac{R^{TM}_t}{4\pi}+61\left(\frac{1}{2R_t^{TM}}-\frac{1}{8\pi{\rm
tan}{\frac{R^{TM}_t}{4\pi}}}\right)\right)\right]dt\right\}^{(11)}.\eqno(3.22)$$
Plugging (3.21) and (3.22) into (3.20) and comparing the constant
terms of both sides, we obtain that
$$\left\{\int_0^1\frac{\sqrt{2}\widehat{L}(TM,\nabla^{TM}_t)}{{\rm
cosh}^2{\frac{c}{2}}}{\rm
tr}\left[A\left(\frac{1}{2R_t^{TM}}-\frac{1}{4\pi{\rm
sin}{\frac{R^{TM}_t}{2\pi}}}\right)\right]\right\}^{(11)}=2^3(2^6z_0+z_1),$$
so we have the following $11$-dimensional analogue of the twisted miraculous cancellation formula.\\

\noindent {\bf Corollary 3.3} {\it The following equality holds}
$$\left\{\int_0^1\frac{\sqrt{2}\widehat{L}(TM,\nabla^{TM}_t)}{{\rm
cosh}^2{\frac{c}{2}}}{\rm
tr}\left[A\left(\frac{1}{2R_t^{TM}}-\frac{1}{4\pi{\rm
sin}{\frac{R^{TM}_t}{2\pi}}}\right)\right]\right\}^{(11)}$$
$$=8\left\{\int_0^1\widehat{A}(TM,\nabla^{TM}_t){\rm
cosh}(\frac{c}{2})\left({\rm
ch}(T_CM,\nabla^{T_CM}_t)-3(e^c+e^{-c}-2)\right)\right.$$
$$\left.\times{\rm tr}\left[A\left(\frac{1}{2R_t^{TM}}-\frac{1}{8\pi{\rm
tan}{\frac{R^{TM}_t}{4\pi}}}\right)\right]dt\right.+$$
$$\left.\int_0^1\widehat{A}(TM,\nabla^{TM}_t){\rm
cosh}(\frac{c}{2}){\rm tr}\left[A\left(-\frac{1}{2\pi}{\rm
sin}\frac{R^{TM}_t}{2\pi}-3\left(\frac{1}{2R_t^{TM}}-\frac{1}{8\pi{\rm
tan}{\frac{R^{TM}_t}{4\pi}}}\right)\right)\right]dt\right\}^{(11)}.\eqno(3.23)$$\\

\indent Next we consider the transgression of
$\Phi_L(\nabla^{TM},\nabla^{\xi},\tau),~\Phi_W(\nabla^{TM},\nabla^{\xi},\tau)$,\\~
\noindent $\Phi_W'(\nabla^{TM},\nabla^{\xi},\tau)$ about
$\nabla^{\xi}$. Let $\nabla_1^{\xi},~\nabla_0^{\xi}$ be two
Euclidean connections on $\xi$ and
$B=\nabla_1^{\xi}-\nabla_0^{\xi}$. By (3.6)-(3.9), we
have\\
 $\Phi_L(\nabla^{TM},\nabla_1^{\xi},\tau)
-\Phi_L(\nabla^{TM},\nabla_0^{\xi},\tau)$ $$=
\frac{1}{8\pi^2}d\int_0^1
\Phi_L(\nabla^{TM},\nabla_t^{\xi},\tau){\rm tr}\left[B \left(
\frac{\theta'_2(\frac{R_t^{\xi}}{4\pi^2},\tau)}
{\theta_2(\frac{R_t^{\xi}}{4\pi^2},\tau)}+\frac{\theta_3'(\frac{R_t^{\xi}}{4\pi^2},\tau)}
{\theta_3(\frac{R_t^{\xi}}{4\pi^2},\tau)}-2\frac{\theta_1'(\frac{R_t^{\xi}}{4\pi^2},\tau)}
{\theta_1(\frac{R_t^{\xi}}{4\pi^2},\tau)}\right)\right]dt.\eqno(3.24)$$
We define\\
$CS\Phi_L(\nabla^{TM},\nabla_0^{\xi},\nabla_1^{\xi},\tau)$
$$:= \frac{\sqrt{2}}{8\pi^2}\int_0^1
\Phi_L(\nabla^{TM},\nabla_t^{\xi},\tau){\rm tr}\left[B \left(
\frac{\theta'_2(\frac{R_t^{\xi}}{4\pi^2},\tau)}
{\theta_2(\frac{R_t^{\xi}}{4\pi^2},\tau)}+\frac{\theta_3'(\frac{R_t^{\xi}}{4\pi^2},\tau)}
{\theta_3(\frac{R_t^{\xi}}{4\pi^2},\tau)}-2\frac{\theta_1'(\frac{R_t^{\xi}}{4\pi^2},\tau)}
{\theta_1(\frac{R_t^{\xi}}{4\pi^2},\tau)}\right)\right]dt.\eqno(3.25)$$
which is in $\Omega^{\rm odd}(M,{\bf C})[[q^{\frac{1}{2}}]].$ Since
$M$ is $4k-1$ dimensional,
$\{CS\Phi_L(\nabla^{TM},\nabla_0^{\xi},\nabla_1^{\xi},\tau)\}^{(4k-1)}$
represents an element in $H^{4k-1}(M,{\bf C})[[q^{\frac{1}{2}}]]$.
Similarly, we can compute the transgressed forms for
$\Phi_W,~\Phi_W'$ respectively and define\\
$CS\Phi_W(\nabla^{TM},\nabla_0^{\xi},\nabla_1^{\xi},\tau)$
$$:= \frac{1}{8\pi^2}\int_0^1
\Phi_W(\nabla^{TM},\nabla_t^{\xi},\tau){\rm tr}\left[B \left(
\frac{\theta'_3(\frac{R_t^{\xi}}{4\pi^2},\tau)}
{\theta_3(\frac{R_t^{\xi}}{4\pi^2},\tau)}+\frac{\theta_1'(\frac{R_t^{\xi}}{4\pi^2},\tau)}
{\theta_1(\frac{R_t^{\xi}}{4\pi^2},\tau)}-2\frac{\theta_2'(\frac{R_t^{\xi}}{4\pi^2},\tau)}
{\theta_2(\frac{R_t^{\xi}}{4\pi^2},\tau)}\right)\right]dt,\eqno(3.26)$$
$CS\Phi_W'(\nabla^{TM},\nabla_0^{\xi},\nabla_1^{\xi},\tau)$
$$:=
\frac{1}{8\pi^2}\int_0^1
\Phi_W'(\nabla^{TM},\nabla_t^{\xi},\tau){\rm tr}\left[B \left(
\frac{\theta'_2(\frac{R_t^{\xi}}{4\pi^2},\tau)}
{\theta_2(\frac{R_t^{\xi}}{4\pi^2},\tau)}+\frac{\theta_1'(\frac{R_t^{\xi}}{4\pi^2},\tau)}
{\theta_1(\frac{R_t^{\xi}}{4\pi^2},\tau)}-2\frac{\theta_3'(\frac{R_t^{\xi}}{4\pi^2},\tau)}
{\theta_3(\frac{R_t^{\xi}}{4\pi^2},\tau)}\right)\right]dt,\eqno(3.26)$$
which also lie in $\Omega^{\rm odd}(M,{\bf C})[[q^{\frac{1}{2}}]]$
and their top components represent elements in $H^{4k-1}(M,{\bf
C})[[q^{\frac{1}{2}}]]$. Similarly we have\\

 \noindent {\bf Theorem 3.4} {\it Let $M$ be a $4k-1$ dimensional manifold and
$\nabla^{TM}$ be a connection on $TM$ and $\xi$ be a two dimensional
oriented Euclidean real vector
bundle with two Euclidean connections $\nabla_1^{\xi}$,~$\nabla_0^{\xi}$, then we have\\
\noindent 1)
$\{CS\Phi_L(\nabla^{TM},\nabla_0^{\xi},\nabla_1^{\xi},\tau)\}^{(4k-1)}$
is a modular form of weight $2k$ over $\Gamma_0(2)$;\\
$\{CS\Phi_W(\nabla^{TM},\nabla_0^{\xi},\nabla_1^{\xi},\tau)\}^{(4k-1)}$
is a modular form of weight $2k$ over $\Gamma^0(2);$\\
$\{CS\Phi_W'(\nabla^{TM},\nabla_0^{\xi},\nabla_1^{\xi},\tau)\}^{(4k-1)}$
is a modular form of weight $2k$ over $\Gamma_\theta(2).$\\
2) The following equalities hold,}
$$\{CS\Phi_L(\nabla^{TM},\nabla_0^{\xi},\nabla_1^{\xi},\tau)\}^{(4k-1)}
=(2\tau)^{2k}\{CS\Phi_W(\nabla^{TM},\nabla_0^{\xi},\nabla_1^{\xi},\tau)\}^{(4k-1)},$$
$$CS\Phi_W(\nabla^{TM},\nabla_0^{\xi},\nabla_1^{\xi},\tau)
=CS\Phi_W'(\nabla^{TM},\nabla_0^{\xi},\nabla_1^{\xi},\tau).$$
\noindent{\bf Proof.} By (3.14),(3.16) and
$$\frac{\theta_2'(z,-\frac{1}{\tau})} {\theta_2(z,-\frac{1}{\tau})}
+\frac{\theta_3'(z,-\frac{1}{\tau})} {\theta_3(z,-\frac{1}{\tau})}
-2\frac{\theta_1'(z,-\frac{1}{\tau})} {\theta_1(z,-\frac{1}{\tau})}
=\tau\left(\frac{\theta_1'(\tau z,{\tau})} {\theta_1(\tau z,{\tau})}
+\frac{\theta_3'(\tau z,{\tau})} {\theta_3(\tau z,{\tau})}
-2\frac{\theta_2'(\tau z,{\tau})} {\theta_2(\tau
z,{\tau})}\right),\eqno(3.27)$$
we can get
$$\{CS\Phi_L(\nabla^{TM},\nabla_0^{\xi},\nabla_1^{\xi},-\frac{1}{\tau})\}^{(4k-1)}
=(2\tau)^{2k}\{CS\Phi_W(\nabla^{TM},\nabla_0^{TM},\nabla_1^{\xi},\tau)\}^{(4k-1)},\eqno(3.28)$$
Similarly we can show that
$$CS\Phi_L(\nabla^{TM},\nabla_0^{\xi},\nabla_1^{\xi},{\tau}+1)
=CS\Phi_L(\nabla^{TM},\nabla_0^{\xi},\nabla_1^{\xi},\tau),$$
$$\{CS\Phi_W(\nabla^{TM},\nabla_0^{\xi},\nabla_1^{\xi},-\frac{1}{\tau})\}^{(4k-1)}
=(\frac{\tau}{2})^{2k}\{CS\Phi_L(\nabla^{TM},\nabla_0^{\xi},\nabla_1^{\xi},\tau)\}^{(4k-1)},$$
$$CS\Phi_W(\nabla^{TM},\nabla_0^{\xi},\nabla_1^{\xi},{\tau}+1)
=CS\Phi_W'(\nabla^{TM},\nabla_0^{\xi},\nabla_1^{\xi},\tau),$$
$$\{CS\Phi_W'(\nabla^{TM},\nabla_0^{\xi},\nabla_1^{\xi},-\frac{1}{\tau})\}^{(4k-1)}
=(\tau)^{2k}\{CS\Phi_W'(\nabla^{TM},\nabla_0^{\xi},\nabla_1^{\xi},\tau)\}^{(4k-1)},$$
$$CS\Phi_W'(\nabla^{TM},\nabla_0^{\xi},\nabla_1^{\xi},{\tau}+1)
=CS\Phi_W(\nabla^{TM},\nabla_0^{\xi},\nabla_1^{\xi},\tau).\eqno(3.29)$$
From (3.28) and (3.29), we can get
$\{CS\Phi_L(\nabla^{TM},\nabla_0^{\xi},\nabla_1^{\xi},{\tau})\}^{(4k-1)}$
is a modular form of weight $2k$ over $\Gamma_0(2).$ Similarly we
can prove that
$\{CS\Phi_W(\nabla^{TM},\nabla_0^{\xi},\nabla_1^{\xi},\tau)\}^{(4k-1)}$
is a modular form of weight $2k$ over $\Gamma^0(2)$ and
$\{CS\Phi_W'(\nabla^{TM},\nabla_0^{\xi},\nabla_1^{\xi},\tau)\}^{(4k-1)}$
is a modular form of weight $2k$ over $\Gamma_\theta(2).$ $\Box$\\
\indent Let $M$ be a compact oriented smooth $3$-dimensional
manifold, we have
\begin{eqnarray*}
&~&CS\Phi_L(\nabla^{TM},\nabla_0^{\xi},\nabla_1^{\xi},\tau)\\
&=& \frac{\sqrt{2}}{8\pi^2}\int_0^1
\Phi_L(\nabla^{TM},\nabla_t^{\xi},\tau){\rm tr}\left[B \left(
\frac{\theta'_2(\frac{R_t^{\xi}}{4\pi^2},\tau)}
{\theta_2(\frac{R_t^{\xi}}{4\pi^2},\tau)}+\frac{\theta_3'(\frac{R_t^{\xi}}{4\pi^2},\tau)}
{\theta_3(\frac{R_t^{\xi}}{4\pi^2},\tau)}-2\frac{\theta_1'(\frac{R_t^{\xi}}{4\pi^2},\tau)}
{\theta_1(\frac{R_t^{\xi}}{4\pi^2},\tau)}\right)\right]dt\\
&=& \frac{1}{2\pi^2}\int_0^1 {\rm tr}\left[B \left(
\frac{\theta'_2(\frac{R_t^{\xi}}{4\pi^2},\tau)}
{\theta_2(\frac{R_t^{\xi}}{4\pi^2},\tau)}+\frac{\theta_3'(\frac{R_t^{\xi}}{4\pi^2},\tau)}
{\theta_3(\frac{R_t^{\xi}}{4\pi^2},\tau)}-2\frac{\theta_1'(\frac{R_t^{\xi}}{4\pi^2},\tau)}
{\theta_1(\frac{R_t^{\xi}}{4\pi^2},\tau)}\right)\right]dt\\
&=& \frac{1}{8\pi^4}\frac{\partial}{\partial z}\left(
\frac{\theta'_2(z,\tau)} {\theta_2(z,\tau)}+\frac{\theta_3'(z,\tau)}
{\theta_3(z,\tau)}-2\frac{\theta_1'(z,\tau)}
{\theta_1(z,\tau)}\right)|_{z=0}\int_0^1 {\rm tr}[BR_t^{\xi}]dt.
\end{eqnarray*}
Since $\frac{\partial}{\partial z}\left( \frac{\theta'_2(z,\tau)}
{\theta_2(z,\tau)}+\frac{\theta_3'(z,\tau)}
{\theta_3(z,\tau)}-2\frac{\theta_1'(z,\tau)}
{\theta_1(z,\tau)}\right)|_{z=0}$ is a modular form of weight $2$
over $\Gamma_0(2),$ then it is a scalar multiple of
$\delta_1(\tau)$. Direct computations show
$$\frac{\partial}{\partial z}\left( \frac{\theta'_2(z,\tau)}
{\theta_2(z,\tau)}+\frac{\theta_3'(z,\tau)}
{\theta_3(z,\tau)}-2\frac{\theta_1'(z,\tau)}
{\theta_1(z,\tau)}\right)|_{z=0}=2\pi^2+O(q^{\frac{1}{2}}),$$ so
$$\frac{\partial}{\partial z}\left(
\frac{\theta'_2(z,\tau)} {\theta_2(z,\tau)}+\frac{\theta_3'(z,\tau)}
{\theta_3(z,\tau)}-2\frac{\theta_1'(z,\tau)}
{\theta_1(z,\tau)}\right)|_{z=0}=8\pi^2\delta_1(\tau).$$ By (4.15)
in [CH], we have
$$CS\Phi_L(\nabla^{TM},\nabla_0^{\xi},\nabla_1^{\xi},\tau)
=\frac {1}{2\pi^2}\delta_1(\tau){\rm
tr}\left[B[\nabla_0^{\xi},\nabla_1^{\xi}]+\frac{2}{3}B\wedge B\wedge
B\right].\eqno(3.30)$$ Similarly, we obtain that
$$CS\Phi_W(\nabla^{TM},\nabla_0^{\xi},\nabla_1^{\xi},\tau)
=\frac {1}{8\pi^2}\delta_2(\tau){\rm
tr}\left[B[\nabla_0^{\xi},\nabla_1^{\xi}]+\frac{2}{3}B\wedge B\wedge
B\right],\eqno(3.31)$$
$$CS\Phi'_W(\nabla^{TM},\nabla_0^{\xi},\nabla_1^{\xi},\tau)
=\frac {1}{8\pi^2}\delta_3(\tau){\rm
tr}\left[B[\nabla_0^{\xi},\nabla_1^{\xi}]+\frac{2}{3}B\wedge B\wedge
B\right].\eqno(3.32)$$
 \indent Let $M$ be $11$ dimensional and
$k=3$. We have that
$\{CS\Phi_L(\nabla^{TM},\nabla_0^{\xi},\nabla_1^{\xi},\tau)\}^{(11)}$
is a modular form of weight $6$ over $\Gamma_0(2),$
$\{CS\Phi_W(\nabla^{TM},\nabla_0^{\xi},\nabla_1^{\xi},\tau)\}^{(11)}$
is a modular form of weight $6$ over $\Gamma^0(2)$ and
$$\{CS\Phi_L(\nabla^{TM},\nabla_0^{\xi},\nabla_1^{\xi},-\frac{1}{\tau})\}^{(11)}
=(2\tau)^{6}\{CS\Phi_W(\nabla^{TM},\nabla_0^{\xi},\nabla_1^{\xi},\tau)\}^{(11)}.$$
By Lemma 2.2, we have
$$\{CS\Phi_W(\nabla^{TM},\nabla_0^{\xi},\nabla_1^{\xi},\tau)\}^{(11)}
=z_0(8\delta_2)^3+z_1(8\delta_2)\varepsilon_2,\eqno(3.33)$$ and by
(2.19) and Theorem 3.4,
$$\{CS\Phi_L(\nabla^{TM},\nabla_0^{\xi},\nabla_1^{\xi},\tau)\}^{(11)}
=2^6[z_0(8\delta_1)^3+z_1(8\delta_1)\varepsilon_1].\eqno(3.34)$$ By
comparing the $q^{\frac{1}{2}}$-expansion coefficients in (3.33), we
get
$$z_0=\left\{\int_0^1\widehat{A}(TM,\nabla^{TM}){\rm
cos}(\frac{R_t^\xi}{4\pi}){\rm tr}\left[\frac{B}{8\pi}{\rm
tan}\frac{R_t^\xi}{4\pi}\right]dt\right\}^{(11)},\eqno(3.35)$$

$$z_1=\left\{\int_0^1\widehat{A}(TM,\nabla^{TM}){\rm
cos}(\frac{R_t^\xi}{4\pi })\left(3{\rm
ch}(\xi_C,\nabla^{\xi_C}_t)-{\rm ch}(T_CM,\nabla^{T_CM})+77\right)
\right.$$
$$\left.\times {\rm tr}\left[\frac{B}{8\pi}{\rm
tan}\frac{R_t^\xi}{4\pi}\right]dt+\int_0^1\widehat{A}(TM,\nabla^{TM}){\rm
cos}(\frac{R_t^\xi}{4\pi }){\rm tr}\left[\frac{3B}{2\pi}{\rm
sin}\frac{R_t^\xi}{2\pi}\right]dt \right\}^{(11)}.\eqno(3.36)$$
Plugging (3.35) and (3.36) into (3.34) and comparing the constant
terms of both sides, we obtain that\\

\noindent {\bf Corollary 3.5} {\it The following equality holds}
$$\left\{\int_0^1\frac{\widehat{L}(TM,\nabla^{TM})}{{\rm
cos}^2(\frac{R_t^\xi}{4\pi})}{\rm tr}\left[B{\rm
tan}\frac{R_t^\xi}{4\pi}\right]dt\right\}^{(11)} $$
$$=16\sqrt{2}\pi
\left\{\int_0^1\widehat{A}(TM,\nabla^{TM}){\rm
cos}(\frac{R_t^\xi}{4\pi })\left(3{\rm
ch}(\xi_C,\nabla^{\xi_C}_t)-{\rm ch}(T_CM,\nabla^{T_CM})+13\right)
\right.$$
$$\left.\times {\rm tr}\left[\frac{B}{8\pi}{\rm
tan}\frac{R_t^\xi}{4\pi}\right]dt+\int_0^1\widehat{A}(TM,\nabla^{TM}){\rm
cos}(\frac{R_t^\xi}{4\pi }){\rm tr}\left[\frac{3B}{2\pi}{\rm
sin}\frac{R_t^\xi}{2\pi}\right]dt \right\}^{(11)}.\eqno(3.37)$$

\section {Transgressed forms and modularities on $4k+1$ dimensional
manifolds}

   \quad Let $M$ be a $4k+1$ dimensional Riemannian manifold. Set
   $$\Theta_1(T_{C}M+\xi_C,\xi_C)=\bigotimes _{n=1}^{\infty}S_{q^n}(\widetilde{T_CM+\xi_C})\otimes
\bigotimes
_{m=1}^{\infty}\wedge_{q^m}(\widetilde{T_CM+\xi_C}-2\widetilde{\xi_C})$$
$$~~~~~~~\otimes \bigotimes _{r=1}^{\infty}\wedge
_{q^{r-\frac{1}{2}}}(\widetilde{\xi_C})\otimes\bigotimes
_{s=1}^{\infty}\wedge _{-q^{s-\frac{1}{2}}}(\widetilde{\xi_C}),$$
$$\Theta_2(T_{C}M+\xi_C,\xi_C)=\bigotimes _{n=1}^{\infty}S_{q^n}(\widetilde{T_CM+\xi_C})\otimes
\bigotimes
_{m=1}^{\infty}\wedge_{-q^{m-\frac{1}{2}}}(\widetilde{T_CM+\xi_C}-2\widetilde{\xi_C})$$
$$~~~~~~\otimes \bigotimes _{r=1}^{\infty}\wedge
_{q^{r-\frac{1}{2}}}(\widetilde{\xi_C})\otimes\bigotimes
_{s=1}^{\infty}\wedge _{q^{s}}(\widetilde{\xi_C}),$$
$$\Theta_3(T_{C}M+\xi_C,\xi_C)=\bigotimes _{n=1}^{\infty}S_{q^n}(\widetilde{T_CM+\xi_C})\otimes
\bigotimes
_{m=1}^{\infty}\wedge_{q^{m-\frac{1}{2}}}(\widetilde{T_CM+\xi_C}-2\widetilde{\xi_C})$$
$$~~~~~~\otimes \bigotimes _{r=1}^{\infty}\wedge
_{q^{r}}(\widetilde{\xi_C})\otimes\bigotimes _{s=1}^{\infty}\wedge
_{-q^{s-\frac{1}{2}}}(\widetilde{\xi_C}).\eqno(4.1)$$ Define
\begin{eqnarray*}
\widetilde{\Phi_L}(\nabla^{TM},\nabla^{\xi},\tau)&=&\widehat{L}(TM,\nabla^{TM})\frac{{\rm
cosh}(\frac{c}{2})}{{\rm sinh}(\frac{c}{2})}\\
&~&\cdot\left({\rm ch}(\Theta_1(T_{C}M+\xi_C,C^2))-\frac{{\rm
ch}(\Theta_1(T_{C}M+\xi_C,\xi_C))}{{\rm
cosh}^2(\frac{c}{2})}\right),\\
\widetilde{\Phi_W}(\nabla^{TM},\nabla^{\xi},\tau)&=&{\widehat{A}(TM,\nabla^{TM})}\frac{1}{2{\rm
sinh}(\frac{c}{2})}\left({\rm ch}(\Theta_2(T_{C}M+\xi_C,C^2))\right.\\
&~&\left.-{\rm cosh}(\frac{c}{2}){\rm ch}(\Theta_2(T_{C}M+\xi_C,\xi_C))\right)\\
\widetilde{\Phi'_W}(\nabla^{TM},\nabla^{\xi},\tau)&=&{\widehat{A}(TM,\nabla^{TM})}\frac{1}{2{\rm
sinh}(\frac{c}{2})}\left({\rm ch}(\Theta_3(T_{C}M+\xi_C,C^2))\right.\\
&~&\left.-{\rm cosh}(\frac{c}{2}){\rm
ch}(\Theta_3(T_{C}M+\xi_C,\xi_C))\right)~~~~~~~~~~~~~~~~~~~~~~(4.2)
\end{eqnarray*}
 Through direct computations, we
get (cf. [HH])
\begin{eqnarray*}
\widetilde{\Phi_L}(\nabla^{TM},\nabla^{\xi},\tau)&=&\frac{\sqrt{2}^{4k+1}}{\pi\sqrt{-1}}\left(
\prod_{j=1}^{2k}x_j\frac{\theta'(0,\tau)}{\theta(x_j,\tau)}\frac
{\theta_1(x_j,\tau)}{\theta_1(0,\tau)}\right)
\frac{\theta'(0,\tau)}{\theta(u,\tau)}\\
&~&\cdot\left(\frac {\theta_1(u,\tau)}{\theta_1(0,\tau)} -\frac
{\theta_1(0,\tau)}{\theta_1(u,\tau)} \frac
{\theta_3(u,\tau)}{\theta_3(0,\tau)} \frac
{\theta_2(u,\tau)}{\theta_2(0,\tau)}\right)\\
\widetilde{\Phi_W}(\nabla^{TM},\nabla^{\xi},\tau)&=&\frac{1}{2\pi\sqrt{-1}}\left(
\prod_{j=1}^{2k}x_j\frac{\theta'(0,\tau)}{\theta(x_j,\tau)}\frac
{\theta_2(x_j,\tau)}{\theta_2(0,\tau)}\right)
\frac{\theta'(0,\tau)}{\theta(u,\tau)}\\
&~&\cdot\left(\frac {\theta_2(u,\tau)}{\theta_2(0,\tau)} -\frac
{\theta_2(0,\tau)}{\theta_2(u,\tau)} \frac
{\theta_3(u,\tau)}{\theta_3(0,\tau)} \frac
{\theta_1(u,\tau)}{\theta_1(0,\tau)}\right)\\
\widetilde{\Phi_W}(\nabla^{TM},\nabla^{\xi},\tau)&=&\frac{1}{2\pi\sqrt{-1}}\left(
\prod_{j=1}^{2k}x_j\frac{\theta'(0,\tau)}{\theta(x_j,\tau)}\frac
{\theta_3(x_j,\tau)}{\theta_3(0,\tau)}\right)
\frac{\theta'(0,\tau)}{\theta(u,\tau)}\\
&~&\cdot\left(\frac {\theta_3(u,\tau)}{\theta_3(0,\tau)} -\frac
{\theta_3(0,\tau)}{\theta_3(u,\tau)} \frac
{\theta_2(u,\tau)}{\theta_2(0,\tau)} \frac
{\theta_1(u,\tau)}{\theta_1(0,\tau)}\right)~~~~~~(4.3)
\end{eqnarray*}
 Applying the Chern-Weil
theory and Lemma 3.1 again, we can transgress
$\widetilde{\Phi_L},~\widetilde{\Phi_W}, \widetilde{\Phi_W'}$ about
$\nabla_{TM}$ and define transgressed forms as follows:\\
$CS\widetilde{\Phi_L}(\nabla_0^{TM},\nabla_1^{TM},\nabla^{\xi},\tau)$
$$:= \frac{\sqrt{2}}{8\pi^2}\int_0^1
\widetilde{\Phi_L}(\nabla_t^{TM},\nabla^{\xi},\tau){\rm tr}\left[A
\left(
\frac{1}{\frac{R_t^{TM}}{4\pi^2}}-\frac{\theta'(\frac{R_t^{TM}}{4\pi^2},\tau)}
{\theta(\frac{R_t^{TM}}{4\pi^2},\tau)}+\frac{\theta_1'(\frac{R_t^{TM}}{4\pi^2},\tau)}
{\theta_1(\frac{R_t^{TM}}{4\pi^2},\tau)}\right)\right]dt.\eqno(4.4)$$
which is in $\Omega^{\rm odd}(M,{\bf C})[[q^{\frac{1}{2}}]].$ Since
$M$ is $4k+1$ dimensional,
$\{CS\widetilde{\Phi_L}(\nabla_0^{TM},\nabla_1^{TM},\nabla^{\xi},\tau)\}^{(4k+1)}$
represents an element in $H^{4k+1}(M,{\bf C})[[q^{\frac{1}{2}}]]$.
Similarly, we can define\\
$CS\widetilde{\Phi_W}(\nabla_0^{TM},\nabla_1^{TM},\nabla^{\xi},\tau)$
$$:= \frac{1}{8\pi^2}\int_0^1
\widetilde{\Phi_W}(\nabla_t^{TM},\nabla^{\xi},\tau){\rm tr}\left[A
\left(
\frac{1}{\frac{R_t^{TM}}{4\pi^2}}-\frac{\theta'(\frac{R_t^{TM}}{4\pi^2},\tau)}
{\theta(\frac{R_t^{TM}}{4\pi^2},\tau)}+\frac{\theta_2'(\frac{R_t^{TM}}{4\pi^2},\tau)}
{\theta_2(\frac{R_t^{TM}}{4\pi^2},\tau)}\right)\right]dt;\eqno(4.5)$$
$CS\widetilde{\Phi_W'}(\nabla_0^{TM},\nabla_1^{TM},\nabla^{\xi},\tau)$
$$:= \frac{1}{8\pi^2}\int_0^1
\widetilde{\Phi_W'}(\nabla_t^{TM},\nabla^{\xi},\tau){\rm tr}\left[A
\left(
\frac{1}{\frac{R_t^{TM}}{4\pi^2}}-\frac{\theta'(\frac{R_t^{TM}}{4\pi^2},\tau)}
{\theta(\frac{R_t^{TM}}{4\pi^2},\tau)}+\frac{\theta_3'(\frac{R_t^{TM}}{4\pi^2},\tau)}
{\theta_3(\frac{R_t^{TM}}{4\pi^2},\tau)}\right)\right]dt.\eqno(4.6)$$
Using the same discussions as Theorem 3.2, we obtain\\

 \noindent
{\bf Theorem 4.1} {\it Let $M$ be a $4k+1$ dimensional manifold and
$\nabla_0^{TM},~\nabla_1^{TM}$ be two connections on $TM$ and $\xi$
be a two dimensional oriented Euclidean real vector
bundle with a Euclidean connection $\nabla^{\xi}$, then we have\\
\noindent 1)
$\{CS\widetilde{\Phi_L}(\nabla_0^{TM},\nabla_1^{TM},\nabla^{\xi},\tau)\}^{(4k+1)}$
is a modular form of weight $2k+2$ over $\Gamma_0(2)$;\\
$\{CS\widetilde{\Phi_W}(\nabla_0^{TM},\nabla_1^{TM},\nabla^{\xi},\tau)\}^{(4k+1)}$
is a modular form of weight $2k+2$ over $\Gamma^0(2);$\\
$\{CS\widetilde{\Phi_W'}(\nabla_0^{TM},\nabla_1^{TM},\nabla^{\xi},\tau)\}^{(4k+1)}$
is a modular form of weight $2k+2$ over $\Gamma_\theta(2).$\\
2) The following equalities hold,}
$$\{CS\widetilde{\Phi_L}(\nabla_0^{TM},\nabla_1^{TM},\nabla^{\xi},-\frac{1}{\tau})\}^{(4k+1)}
=(2\tau)^{2k+2}\{CS\widetilde{\Phi_W}(\nabla_0^{TM},\nabla_1^{TM},\nabla^{\xi},\tau)\}^{(4k+1)},$$
$$CS\widetilde{\Phi_W}(\nabla_0^{TM},\nabla_1^{TM},\nabla^{\xi},\tau+1)
=CS\widetilde{\Phi_W'}(\nabla_0^{TM},\nabla_1^{TM},\nabla^{\xi},\tau).$$

Let $M$ be $9$ dimensional and $k=2$. We have that
$\{CS\widetilde{\Phi_L}(\nabla_0^{TM},\nabla_1^{TM},\nabla^{\xi},\tau)\}^{(9)}$
is a modular form of weight $6$ over $\Gamma_0(2),$
$\{CS\widetilde{\Phi_W}(\nabla_0^{TM},\nabla_1^{TM},\nabla^{\xi},\tau)\}^{(9)}$
is a modular form of weight $6$ over $\Gamma^0(2)$ and
$$\{CS\widetilde{\Phi_L}(\nabla_0^{TM},\nabla_1^{TM},\nabla^{\xi},-\frac{1}{\tau})\}^{(9)}
=(2\tau)^{6}\{CS\widetilde{\Phi_W}(\nabla_0^{TM},\nabla_1^{TM},\nabla^{\xi},\tau)\}^{(9)}.$$
By Lemma 2.2, we have
$$\{CS\widetilde{\Phi_W}(\nabla_0^{TM},\nabla_1^{TM},\nabla^{\xi},\tau)\}^{(9)}
=z_0(8\delta_2)^3+z_1(8\delta_2)\varepsilon_2,\eqno(4.7)$$ and by
(2.19) and Theorem 4.1,
$$\{CS\widetilde{\Phi_L}(\nabla_0^{TM},\nabla_1^{TM},\nabla^{\xi},\tau)\}^{(9)}
=2^6[z_0(8\delta_1)^3+z_1(8\delta_1)\varepsilon_1].\eqno(4.8)$$ By
comparing the $q^{\frac{1}{2}}$-expansion coefficients in (4.7), we
get
$$z_0=-\left\{\int_0^1\frac{\widehat{A}(TM,\nabla^{TM}_t)}{2{\rm sinh}(\frac{c}{2})}
(1-{\rm cosh}\frac{c}{2}){\rm
tr}\left[A\left(\frac{1}{2R_t^{TM}}-\frac{1}{8\pi{\rm
tan}{\frac{R^{TM}_t}{4\pi}}}\right)\right]dt\right\}^{(9)},\eqno(4.9)$$
$$z_1=\left\{-\int_0^1\frac{\widehat{A}(TM,\nabla^{TM}_t)}{2{\rm sinh}(\frac{c}{2})}
(1-{\rm cosh}\frac{c}{2}){\rm tr}\left[\frac{A}{2\pi}{\rm
sin}\frac{R^{TM}_t}{2\pi}\right]dt\right.$$
$$+\int_0^1\frac{\widehat{A}(TM,\nabla^{TM}_t)}{2{\rm sinh}(\frac{c}{2})}
{\rm tr}\left[A\left(\frac{1}{2R_t^{TM}}-\frac{1}{8\pi{\rm
tan}{\frac{R^{TM}_t}{4\pi}}}\right)\right]$$
$$\left.\cdot\left((1-{\rm cosh}\frac{c}{2})({\rm
ch}(T_CM,\nabla_t^{T_CM})+61)+(1+2{\rm
cosh}\frac{c}{2})(e^c+e^{-c}-2)\right)dt\right\}^{(9)}.\eqno(4.10)$$
Plugging (4.9) and (4.10) into (4.8) and comparing the constant
terms of both sides, we obtain that\\

 \noindent {\bf Corollary 4.2}
{\it The following equality holds}
$$\left\{\int_0^1
\sqrt{2}\widehat{L}(TM,\nabla^{TM}_t)\frac{{\rm sinh}\frac{c}{2}}
{{\rm cosh}{\frac{c}{2}}}{\rm
tr}\left[A\left(\frac{1}{2R_t^{TM}}-\frac{1}{4\pi{\rm
sin}{\frac{R^{TM}_t}{2\pi}}}\right)\right]\right\}^{(9)}$$
$$=8\left\{-\int_0^1\frac{\widehat{A}(TM,\nabla^{TM}_t)}{2{\rm sinh}(\frac{c}{2})}
(1-{\rm cosh}\frac{c}{2}){\rm tr}\left[\frac{A}{2\pi}{\rm
sin}\frac{R^{TM}_t}{2\pi}\right]dt\right.$$
$$+\int_0^1\frac{\widehat{A}(TM,\nabla^{TM}_t)}{2{\rm sinh}(\frac{c}{2})}
{\rm tr}\left[A\left(\frac{1}{2R_t^{TM}}-\frac{1}{8\pi{\rm
tan}{\frac{R^{TM}_t}{4\pi}}}\right)\right]$$
$$\left.\cdot\left((1-{\rm cosh}\frac{c}{2})({\rm
ch}(T_CM,\nabla_t^{T_CM})-3)+(1+2{\rm
cosh}\frac{c}{2})(e^c+e^{-c}-2)\right)dt\right\}^{(9)}.\eqno(4.11)$$\\

\noindent{\bf Acknowledgement}~This work was supported by NSFC No.
10801027.\\

\noindent {\bf References}\\

\noindent [AW] L. Alvarez-Gaum\'{e}, E. Witten, Graviational
anomalies,
{\it Nucl. Phys.} B234 (1983), 269-330.\\
\noindent [Ch] K. Chandrasekharan, {\it Elliptic Functions},
Spinger-Verlag, 1985. \\
\noindent [CH] Q. Chen, F. Han, Elliptic genera, transgression and
loop space Chern-Simons form, arXiv:0605366.\\
\noindent [HZ1]F. Han, W. Zhang, ${\rm Spin}^c$-manifold and
elliptic genera,
{\it C. R. Acad. Sci. Paris Serie I.,} 336 (2003), 1011-1014.\\
\noindent [HZ2] F. Han, W. Zhang, Modular invariance, characteristic
numbers and eta Invariants,  {\it J.
Diff. Geom.} 67 (2004), 257-288.\\
 \noindent [HH] F. Han, X. Huang, Even dimensional manifolds and generalized anomaly cancellation
formulas , {\it Trans. AMS }359 (2007), No. 11, 5365-5382.\\
\noindent [Li] K. Liu, Modular invariance and characteristic
numbers. {\it Commu.Math. Phys.} 174 (1995), 29-42.\\
\noindent [Z] W. Zhang, {\it Lectures on Chern-weil Theory and
Witten Deformations.} Nankai Tracks in Mathematics Vol. 4, World
Scientific, Singapore, 2001.\\

\indent School of Mathematics and Statistics , Northeast
Normal University, Changchun, Jilin 130024, China ;\\

 \indent E-mail: {\it wangy581@nenu.edu.cn}
\end {document}